%% file: main.tex
\documentclass[smallextended]{svjour3}       %
    \smartqed

    \makeatletter
    \def\cl@chapter{\@elt {theorem}}
    \makeatother

\input{setup.tex}
\journalname{Set-Valued and Variational Analysis}
\begin{document}

\title{Nonsmoothness in Machine Learning: specific structure, proximal identification, and applications}

\titlerunning{Nonsmoothness in Machine Learning}        %

\author{Franck Iutzeler \and
        Jérôme Malick }

\institute{F. Iutzeler \at
              Univ. Grenoble Alpes 
           \and
           J. Malick \at
              CNRS, LJK
}

\date{\today}

\maketitle

\begin{abstract}
Nonsmoothness is often a curse for optimization; but it is sometimes a blessing, in particular for applications in machine learning. In this paper, we present the specific structure of nonsmooth optimization problems appearing in machine learning and illustrate how to leverage this structure in practice, for compression, acceleration, or dimension reduction. 
We pay a special attention to the presentation to make it concise and easily accessible, with both simple examples and general results.
\end{abstract}

\section{Introduction}

Optimization is at the core of machine learning and nonsmoothness plays there a special role. Often, nonsmooth functions are introduced in learning problems to enforce low-complexity of the optimal model parameters. This promoted structure (e.g.\;sparsity, low-rank, or controlled variation) turns out to be progressively identified by proximal algorithms; most interestingly, it can be leveraged numerically to improve optimization methods. This situation contrasts with a large part of the optimization literature where dealing with nonsmooth functions or checking nonsmooth substructure is usually uneasy.

In this paper, we lay out a generic framework of optimization problems for which nonsmoothness can be exploited, covering most machine learning problems with low-complexity priors. We then go over the main technical tools that enable to mathematically and practically handle nonsmooth structures. We pay a special attention to being accessible for a wide audience\footnote{Advanced readers can jump directly to Sections~\ref{sec:ident} (about general identification) and~\ref{sec:harness} (with entry points to recent research on this topic).} in optimization, including graduate students in applied maths.
No specific knowledge on machine learning or nonsmooth analysis is necessary to follow our developments, but fundamentals in mathematical optimization are required. 
Let us mention that we will make a constant (but basic) use of the {proximal operator} 
which is a key tool to deal with explicit nonsmoothness. 
For a function $g\colon\RR^n\rightarrow \RR \cup\{+\infty\}$ and a parameter\;$\gamma>0$, we define $\prox_{\gamma g}(u)$ for any $u\in\mathbb{R}^n$ as the optimal solution of the following problem: 
\begin{align}\tag{prox}
    \label{eq:prox}
     \displaystyle\prox_{\gamma g} (u) := \arg\min_{y\in\mathbb{R}^n} \left\{ g(y)  + \frac{1}{2\gamma} \| y-u\|_2^2 \right\} .
\end{align}
This properly defines an operator $\prox_{\gamma g}\colon\RR^n\rightarrow \RR^n$ in many cases, in particular when $g$ is convex. We refer to the textbooks\;\cite{hull} and\;\cite{bauschke2011convex} or to the review\;\cite{parikh2014proximal} for more information and an historical perspective.

The rest of the paper is organized as follows. In Section~\ref{sec:learning}, we explain why nonsmooth models are considered in machine learning. In Section~\ref{sec:NS}, we formalize a large class of nonsmooth optimization problems for which structure can be exploited. Then, in Section~\ref{sec:algo}, we discuss the proximal optimization algorithms that can harness this structure. In Section~\ref{sec:ident}, we lay out the mathematical grounds that enable to properly identify structure. Finally, in Section~\ref{sec:harness}, we review recent advances in the use of identified structure to improve performances of some optimization algorithms in machine learning.

\section{Nonsmooth problems in Machine Learning}
\label{sec:learning}

Standard machine learning tasks such as regression, classification, or clustering (see e.g the textbook\;\cite{shalev2014understanding}) lead to optimization problems. 
Indeed, many objectives in supervised learning can write as minimizing some \emph{risk} function $\risk$ measuring the dissimilarity between a learning model, obtained from statistical modeling, with parameters $x\in\mathbb{R}^n$, and a set of $m$ examples $\{a_i,b_i\}_{i=1}^m$, in the form of an input ($a_i$) and target ($b_i$) pair. This problem, often called Empirical Risk Minimization (ERM), has the following form: %
\begin{align}
\tag{ERM}
\label{ERM}
    \min_{x\in\mathbb{R}^n}  ~\risk\left( x ; \{a_i,b_i\}_{i=1}^m \right).
\end{align}

Often, the risk $\risk$ is taken as the average of some loss over the training examples:
$ \risk\left( x ; \{a_i,b_i\}_{i=1}^m \right) = \frac{1}{m}\sum_{i=1}^m \ell( p( x ; a_i ) ; b_i )$
where $p( x ; a )$ is the model prediction for input $a$ using parameters $x$, and $\ell(p;b)$ quantifies the error between the predicted point $p$ and the true target $b$. The simplest case is when the prediction function is linear ($p( x ; a )= \langle  x ; a \rangle$) and the error is quadratic ($\ell(p;b) = (p-b)^2$) leading to the well known least-squares regression problem. Other popular choices for the model include polynomial and gaussian kernels, or neural networks. The error functions often depend either on (i) the statistical modelling of the problem through the (log) likelihood (squared $2$-norm for linear regression of points corrupted by a Gaussian noise, $1$-norm for Laplace noise); or (ii) directly from applications (errors in classification).   %

\begin{figure}[!h]
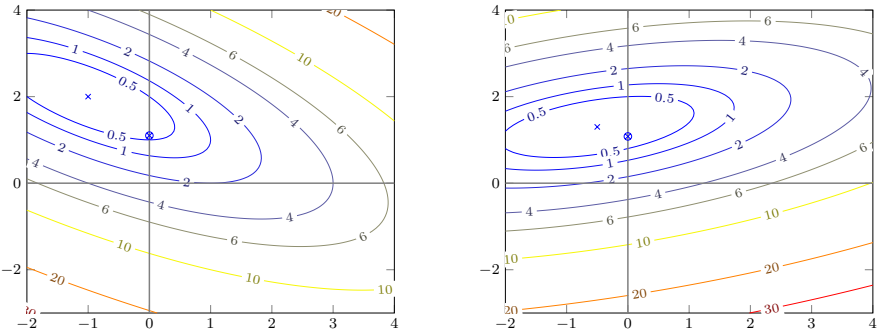

     \centering
     \begin{subfigure}[t]{0.47\textwidth}
         \centering
         \resizebox{\textwidth}{!}{         \input{identif/plotlevel_1.tikz}      }
     \end{subfigure}\hfill
     \begin{subfigure}[t]{0.47\textwidth}
        \centering
        \resizebox{\textwidth}{!}{     \input{identif/plotlevel_2.tikz}   }
    \end{subfigure}
    \caption{Illustration of the stability of optimal solutions of lasso in $\RR^2$. We plot the level sets of $1/2\|Ax-b\|_2^2 $ %
    for two \eqref{lasso} problems with different but close design matrices $A$. We see that while the solutions of the (unregularized) least-squares problem (marked by a x) are different, the solutions of the lasso (marked by a circled cross), although different, lie on the same axis, corresponding to the nonsmoothness loci of the $\ell_1$-norm.\label{fig:lasso}}
\end{figure}

With the ever-growing collection of data, the size of learning problems has significantly increased, both in terms of number of examples, $m$,  as well as in size of optimized parameter, $n$. While the increase in number of examples is generally beneficial for the general conditioning of the problem, the increase of the parameter space often makes the problem ill-conditioned and reduces both the interpretability and stability of the model. In order to overcome this issue, a generally admitted solution is to introduce a \emph{prior} on the \emph{structure} of the model $x$ and to \emph{regularize} the \eqref{ERM} in order to promote the prior structure. One of the first and most well-known instances of such a prior is the sparsity sought in the parameters of least-squares regression, leading to the well known lasso problem \cite{tibshirani1996regression}. %
More generally, regularizing a learning problem to enforce a prior structure conforms to the following steps:
\vspace*{-1ex}
\begin{itemize}
    \item[(i)] \emph{Define a prior structure}. Let us observe the usual priors in machine learning: sparsity, constant by block, fixed rank, etc. They can be described as subsets of $\mathbb{R}^n$ that are rather easy to describe and project on. As mathematical objects, and to comply with a vast part of the literature, we shall see them as (affine or smooth) {manifolds}, but most of the arguments in this paper hold for simple closed sets.

\vspace*{1ex}
    \item[(ii)] \emph{Find an additive nonsmooth function $r$ enforcing this structure}. Nonsmoothness of functions \emph{traps} optimal solutions in low-dimensional manifolds: small perturbations in the risk around these points would not break down optimal structure as illustrated by Figure\;\ref{fig:lasso}. Mathematically, this is due to subdifferentials with non-empty relative interiors in relevant directions. For instance, the subdifferential of the $\ell_1$-norm has a nonempty relative interior along the axes, promoting sparsity. 
\end{itemize}

At this point, a regularized version of \eqref{ERM} can be formulated as
\begin{align}
\tag{Regularized ERM}
\label{reg_ERM}
    \min_{x\in\mathbb{R}^n} ~~ \underbrace{ \risk\left( x ; \{a_i,b_i\}_{i=1}^m \right)}_{\substack{=: f(x) \\[0.5ex] \text{minimizes the risk}} } ~~+ \underbrace{~~ \lambda ~ r(x) ~~ }_{\substack{=: g(x) \\[0.5ex] \text{enforces structure}} } .
\end{align}
To make it more concrete, let us consider the popular example of $\ell_1$-regularized least-squares problems (often called lasso \cite{tibshirani1996regression}): take a linear prediction function, a quadratic loss, %
and the $\ell_1$-norm as a regularizer, then the regularized problem becomes
\begin{align}
\tag{lasso}
\label{lasso}
    \min_{x\in\mathbb{R}^n} ~~ \frac{1}{2}\|A\,x - b \|_2^2 ~+~ \lambda\; \|x\|_1 \;.
\end{align}
Optimal solutions are then sparser than the one of the original least-squares, and their sparsity pattern is stable under small perturbations; see Figure\;\ref{fig:lasso}. 

In general, the regularized risk minimization has thus a  composite form, involving two different parts: $f$ which is linked to the risk itself, and $g$ which only promotes structure, with a hyperparameter $\lambda>0$ controlling the balance between $f$ and $g$. 
This particular formulation enables {one} to use ad-hoc optimization methods which forms the third part of the learning process:
\begin{itemize}
    \item[(iii)] \emph{Design an optimization method complying with this composite formulation.} The structure-enforcing part $g=\lambda \,r$ is necessarily nonsmooth but it may be possible to choose it so that its proximity operator \eqref{eq:prox}
 can be computed easily (with an explicit expression or through a computationally cheap procedure).
When this is the case, proximal methods (i.e. optimization methods involving at least one proximity operator) are the algorithms of choice to solve \eqref{reg_ERM}. Indeed, contrary to subgradient methods, the proximity operator {enables taking} fixed stepsizes and is thus much faster in both theory and practice. We will come back to them in Section~\ref{sec:algo}.
\end{itemize}

To summarize, when solving a machine learning problem, it is common to assume some prior structure and enforce it by adding a regularization term to the \eqref{ERM}. {In the next section, we will see that this prior structure often takes a special and simple form in relation with the proximity operator of $g$.} In this context, we will later explain how proximal methods solving \eqref{reg_ERM} automatically uncover the optimal structure (or at least part of it), and how we can take advantage of this behavior to improve the performance of the algorithms themselves. 

\section{Noticeable Structure in nonsmooth optimization problems}\label{sec:NS}

The type of nonsmoothness encountered in machine learning objectives is often linked to user-defined priors, as in \eqref{reg_ERM}. In this case, the nonsmoothness is under control: it is chosen and relatively simple. 
 In this section, we introduce the notion of Noticeable Structure \eqref{fw} to describe this favorable class of problems for which structure can be {characterized}.

\fwtag{NS}
\begin{fw}[Problems with {Noticeable Structure}]
\label{fw}
For a function $g$ %
and a finite collection\footnote{We mention that the sets in the collection $\collong$ are denoted by $\M_i$ since in most of the literature these sets are described as manifolds. The manifold structure will only be used when needed. However, the closedness of the sets is primordial. Its finiteness is less {essential} but is used to describe properly the structure of a neighborhood points (see the proof of \cref{th:identification}).}
$\collong$ of closed sets, we say that the couple $(g,\col)$ has noticeable structure if 
\vspace*{-2ex}
\begin{itemize}
    \item[1.] \label{fw:1} we have a projection mapping $\proj_{\Mi}$ onto $\Mi$ for all $i$ {as well as for any intersection $\N_{\mathcal{I}} = \cap_{i\in\mathcal{I}} \M_i$ for $\mathcal{I}\subseteq \{1,..,p\}$};
    \item[2.] \label{fw:2} $g$ is non-differentiable at $x$ only if $x\in \Mi$ for some $i$;
    \item[3.] \label{fw:3} $\prox_{\gamma g}(u)$ is a singleton and can be computed explicitly for any $u$ and $\gamma$;
    \item[4.] \label{fw:4} upon computation of $x = \prox_{\gamma g}(u)$, we know if $x\in\Mi$ or not for all $i$.
\end{itemize}
\end{fw}

Before providing examples, let us discuss these four properties. 
Items\;1 and\;2 are important to harness the structure numerically. The first item simply means that it is computationally possible to {project onto any set (or  intersection of sets) present in the collection; this is to enable algorithms to harness the iterates structure by updating preferably along the identified manifold (see Sec.~\ref{sec:alg_ident}).} Item 2 {permits localization} the non-differentiability and as a consequence better grasp differentiability (outside of the manifolds and along manifolds\footnote{In many examples, the nonsmooth function $g$ turns out to be smooth with respect to smooth (sub)manifolds of the collection. We will come back to this relative smoothness of the problem in Section~\ref{sec:riemann}.}).
Item 3 is typically {satisfied} for all proper convex lower semi-continuous functions {and for some classes of non-convex functions (e.g.\;for prox-bounded, prox-regular functions \cite[Prop.\;13.37]{rockafellar-wets-1998}). Note though that (generalized) convexity will not be required for most of our developments.} Item 4 %
is essential: to be able to exploit some structural information, the user needs to \emph{know} the current structure. This condition is not so stringent since for many popular structure-enforcing functions, an explicit proximity operator is known, often based on switch-cases that directly determine the output structure. 

In the remainder of this section, we detail some popular examples of structures and regularizers\footnote{In the following examples, we take $\lambda=1$ out of simplicity since in that case $g\equiv r$. This does not change the discussion on proximity operators since there are explicit formulas for scaling, translation, etc. \cite[Th.~6.1]{beck2017first}. } falling into the framework \eqref{fw} often considered as priors in machine learning \cite{bach2012optimization}.

\subsection{Sparse structure with $\ell_1$-norm}
\label{sec:sparsecol}

The natural manner to look at the sparsity pattern of a vector is to look at each coordinate and see if it is null or not. In terms of  collection of manifolds, we consider 
\begin{align}
\label{eq:sparse}
 \col  = \{\M _1,\ldots,\M _n\} \quad\text{ with }\quad \M _i = \{x\in\mathbb{R}^n  : x_i = 0 \} .
\end{align} 
First, the projections on the $\M_i$'s are direct (\ref{fw}-1). Taking $g(x)= \|x\|_1$, its non-differentiability points exactly match the described manifolds (\ref{fw}-2).

The proximity operator of the $\ell_1$-norm is the well-known {soft-thresholding} operator which can be described coordinate-wise as
\begin{align}
  [x]_i = \left[  \prox_{\gamma \|\cdot\|_1 }(u) \right]_i = \left\{
    \begin{array}{ll}
    ~ 0   & \text{ if } [u]_i \in [ -\gamma , \gamma ]  \\
    ~  [u]_i  - \gamma & \text{ if } [u]_i > \gamma \\
    ~ [u]_i  + \gamma & \text{ if } [u]_i < -\gamma \\
    \end{array}
    \right.
\end{align}
with $[u]_i\in \mathbb{R}$ denoting the $i$-th coordinate of vector\;$u\in \RR^n$. This proximity operator is thus explicitly computable for any $u$ (\ref{fw}-3). Finally, after completing the tests on the right hand side above, the structure of the output $x$ is perfectly known (\ref{fw}-4).

For other types of sparsity-inducing functions and structure (e.g. $\ell_1/\ell_2$ or group sparsity), see \cite{bach2011convex}. Note that these functions do not need to be convex as discussed in the following remark.

\begin{remark}[Sparser structure with non-convex functions] The $\ell_1$-norm is sometimes seen as a convex approximation of the $\ell_0$-``norm''. Although $g(x) = \|x\|_0 = \mathrm{Card}\{ i : [x]_i \neq 0 \}$ is not a norm and is non-convex, it is quite easy to see that it fits the framework \eqref{fw}, with the same collection and the proximal operator
\begin{align}
  [x]_i = \left[  \prox_{\gamma \|\cdot\|_0 }(u) \right]_i = \left\{
    \begin{array}{ll}
    ~ 0   & \text{ if } [u]_i \in [ -\gamma , \gamma ]  \\
    ~  [u]_i  & \text{ if } |[u]_i| > \gamma 
    \end{array}
    \right. ,
\end{align}
sometimes called the {hard-thresholding} operator.

For completeness, notice that all the functions $g(x) = \|x\|_p^p$ for $p\in[0,1]$ {satisfy} \eqref{fw}-1 \& 2 and are thus sparsity-inducing. However, the only explicitly computable proximity operators are for $p=0,0.5,2/3,1$ \cite{chartrand2016nonconvex}, they are thus the only ones matching  \eqref{fw}-3 \& 4. Finally, only $p=1$ (presented above) leads to a convex function.\qed
\end{remark}

\subsection{Flat structure with total variation}
\label{sec:colTV}

The denoising of piecewise-constant one-dimensional\footnote{The 2D problem of matrix regression with flat regions is much harder to define and to solve, see e.g. \cite{condat2017discrete} and references therein.} signals has attracted a lot of attention, notably thanks to the introduction of the total variation function $g(x) = \sum_{i=2}^n |[x]_i - [x]_{i-1}|$ \cite{rudin1992nonlinear}. The structure of piecewise constant signals can be described using the following collection: 
\begin{equation}\label{eq:Mjump}
\col  = \{\M _1,\ldots,\M _{n-1} \} \quad\text{ with }\quad \M _i = \left\{x\in\mathbb{R}^n  :  x_{i} = x_{i-1} \right\},
\end{equation}
and once again \eqref{fw}-1 %
\& 2 %
are easily {satisfied}. The computation of the proximity operator is less direct than before and does not benefit from a closed form solution, however, there are efficient dynamic programming algorithms computing it exactly and enabling to know exactly the structure of the output, see \cite{condat2013direct} for details. Hence, \eqref{fw}-3 \& 4 are also {satisfied}.

\begin{remark}[More structure with non-convex functions] 
In the same vein as before, the total variation can be seen as a convex relaxation of $g(x) = \mathrm{Card} (\{i: [x]_i \neq [x]_{i-1}\})$, sometimes called the {Potts problem}  \cite{weinmann20151}. \eqref{fw}-1 \& 2 are then easily {satisfied}, while 3 \& 4 are less direct but the proximity operator can also be obtained by dynamic programming \cite{friedrich2008complexity}.\qed
\end{remark}

\subsection{Low rank with nuclear norm}
\label{sec:colrank}

One of the most prominent types of structure sought in matrix-values problems is low-rank. Indeed, it is widely used (notably stemming from PCA and spike models) in order to summarize the information in the matrix to its most informative subspaces. It also pops out in matrix factorization and matrix completion problems, with applications in recommender systems. Naturally, the collection writes 
\begin{equation}\label{eq:Mjump}
\col  = \{\M _1,\ldots,\M _{n} \} \quad\text{ with }\quad \M _i = \left\{ \mathrm{rank}(X) = i \right\} 
\end{equation}
and projection onto the manifolds can be obtained numerically by truncating the singular value decomposition (for \ref{fw} 1).

To express the rank of a matrix $X\in\mathbb{R}^{r \times c}$, it is natural to examine its singular values $(\sigma_1,\ldots, \sigma_n)$ (where $n=\min\{r,c\}$). A direct way to promote low-rank matrices is through nuclear norm regularization, which is the $\ell_1$-norm of the singular values vector: $g(X) = \|X\|_* = \sum_{i=1}^n \sigma_i$. One can show that the proximal operator $X = \prox_{\gamma \|\cdot\|_* }(U)$ of a matrix $U = V \diag(\sigma_i) W^*$ can be obtained by 
\begin{align}
 X = V \diag(\nu_i) W^* \text{ with }  \nu_i =   \prox_{\gamma |\cdot| }(\sigma_i)  = \left\{
    \begin{array}{ll}
    ~ 0   & \text{ if } \sigma_i \in [ -\gamma , \gamma ]  \\
    ~  \sigma_i  - \gamma & \text{ if } \sigma_i > \gamma \\
    ~ \sigma_i  + \gamma & \text{ if } \sigma_i < -\gamma \\
    \end{array}
    \right. 
\end{align}
which matches \ref{fw} 3 \& 4.

A remarkable point is that even though two close matrices may have completely different singular vectors, their singular {values} will be close (this is sometimes called Weyl's lemma, see \cite{weyl1912asymptotische,stewart1998perturbation}). Hence, while a new singular value decomposition has to be computed at each application of the proximity operator of the nuclear norm, the rank of the output will be somewhat stable to small perturbations, which is of utmost importance when using sparsity. %

Thus, even though the problem of checking the rank of a matrix can be problematic numerically, the simplicity of the proximity operator (modulo the computation of the SVD) enables it to fall into our framework \eqref{fw}; especially, the rank of the output of the proximity operator is known by construction.

Finally, note that as in the previous case, the ``$\ell_0$-equivalent'' of the nuclear norm is simply the {rank} whose proximity operator involves a hard thresholding of the singular values. Finally, for completeness, several types of matrix structures falling into our framework can be found in \cite{benfenati2018proximal}.

\section{Proximal algorithms}
\label{sec:algo}

From an optimization point of view, regularized empirical risk minimization (see \eqref{reg_ERM} in Section~\ref{sec:learning}) is often seen as a composite problem:%
\begin{equation}
\tag{$\mathcal{P}$}\label{eq:main_problem}
	\min_{x\in \mathbb{R}^n}~f(x)+ g(x)
\end{equation}
where $g$ induces a Noticeable Structure \eqref{fw} as per the framework introduced in the previous section.

The study of optimization algorithms associated with proximity operators, notably in relation with monotone operators \cite{rockafellar1976monotone} and splitting methods \cite{eckstein1992douglas}, has attracted a considerable lot of attention in the optimization community since the 70s; see e.g. \cite{combettes2011proximal,parikh2014proximal} for reviews of the essential points for signal processing and machine learning applications. These proximal algorithms {are written} quite generally as follows with iterations consisting of an $\algo$ step, that defines the algorithm, followed by a proximity operation:
\begin{align}
\tag{$\prox-\algo$}
\label{algo}
\left\{ 
\begin{array}{ll}
    u_{k+1} &= \algo \\
    x_{k+1} &= \prox_{\gamma g}(u_{k+1})
\end{array}
\right.
\end{align}

The main purpose of the $\algo$ step is to handle $f$ which is the part of the problem not managed by $\prox_{\gamma g}$. To do so, it can use the properties of $f$ as well as previous iterates.

We will see in the next section that such proximal methods are the most generic algorithms that can harness nonsmoothness in the framework \eqref{fw}. The only prerequisite is that these methods converge (or converge almost surely, when randomness is involved), which we formalize in the next assumption.

\begin{assumption}
\label{asm:cv}
The method \eqref{algo} is such that
\begin{itemize}
    \item[1.] $(u_k)$ converges (almost surely) to a point $u^\star$,
    \item[2.] $(x_k)$ converges (almost surely) to a point $x^\star$,
    \item[3.] $x^\star = \prox_{\gamma g}(u^\star)$ is a minimizer of \eqref{eq:main_problem}. 
\end{itemize}
\end{assumption}

Let us list some popular proximal algorithms solving \eqref{eq:main_problem} with $g$ convex\footnote{
In the case where $g$ is convex and lower semi-continuous, points \emph{1} and \emph{3} imply point\;\emph{2} and are actually sufficient for convergence. Observe indeed that $\|x_k - x^\star \| = \|\prox_{\gamma g}(x_k) - \prox_{\gamma g}(u^\star) \| \leq \|u_k - u^\star \|$; see more in\;\cite[Chap.~23]{bauschke2011convex}. In non-convex cases, the situation is more complex, all three points are necessary and can be investigated by looking at the prox-regularity and prox-boundedness of $g$ at {a local minimizer} $x^\star$; see \cite{hare2009computing} and references therein.} lower semi-continuous, that {satisfy} \cref{asm:cv}:
\begin{itemize}
    \item \emph{Proximal gradient}: for $f$ convex with an $L$-Lipchitz gradient%
    \begin{align}
    \label{pg}
\tag{PG}
\left\{ 
\begin{array}{ll}
    u_{k+1} &= x_k - \gamma \nabla f(x_k) \\
    x_{k+1} &= \prox_{\gamma g}(u_{k+1})
\end{array}
\right.
\end{align}
   {satisfies}\;\cref{asm:cv} for $\gamma\in(0,2/L)$ {\cite[Th.~10.24]{beck2017first}}.
    
    \vspace*{2ex}
    \item \emph{Accelerated Proximal Gradient}: for $f$ convex with an $L$-Lipchitz gradient,
        \begin{align}
        \label{apg}
\tag{APG}
\left\{ 
\begin{array}{ll}
    u_{k+1} &= x_k + \alpha_k(x_k-x_{k-1}) - \gamma \nabla f( x_k + \alpha_k(x_k-x_{k-1})) \\
    x_{k+1} &= \prox_{\gamma g}(u_{k+1})
\end{array}
\right.
\end{align}
    with $\alpha_k = (k-1)/(k+3)$ {satisfies} \cref{asm:cv} 
    for $\gamma\in(0,1/L]$ {\cite[Th.\;3]{chambolle2015convergence}}.
    \vspace*{2ex}
\item \emph{Douglas-Rachford}: for $f$ convex lower semi-continuous but not necessarily differentiable, the proximity operator of $f$ can be used with a splitting\footnote{A splitting method is necessary here to separate the two proximity operators, since they cannot be directly composed \cite{yu2013decomposing} or averaged \cite{bauschke2008proximal}. } algorithm such as Douglas-Rachford
    \begin{align}
\tag{DR}
\label{dr}
\left\{ 
\begin{array}{ll}
    u_{k+1} &=  \prox_{\gamma f}(2x_k - u_k) + u_k - x_k \\
    x_{k+1} &= \prox_{\gamma g}(u_{k+1})
\end{array}
\right.
\end{align}
which {satisfies}\;\cref{asm:cv}  
for any $\gamma > 0$ \cite[Cor.~27.4]{bauschke2011convex}.

    \vspace*{2ex}
    \item \emph{Variance-Reduced Incremental Methods}: for $f(x) = \frac{1}{m}\sum_{j=1}^m f^j(x)$ with each $(f^j)$ strongly convex with an $L$-Lipchitz gradient,
        \begin{align}
        \label{saga}
\tag{SAGA}
\left\{ 
\begin{array}{ll}
    u_{k+1} &= \displaystyle x_k \!- \!\gamma \Big( \nabla f^{i_k}(x_k) - \nabla f^{i_k}(x_{k-d_k^{i_k}}) + \!\sum_{j=1}^m \nabla f^j\big(x_{k-d^j_k}\big) \!\Big)  \\
    x_{k+1} &= \prox_{\gamma g}(u_{k+1})
\end{array}
\right.
\end{align}
where $i_k$ is drawn uniformly in $\{1,\ldots,m\}$ and $k-d^j_k$ is the last time before $k$ when $\nabla f^{j}$ was computed (i.e. for which $i_{k-d^j_k}=j$). SAGA {satisfies} \cref{asm:cv} for $\gamma=1/(3L)$ {\cite[Cor.~1]{defazio2014saga}}. %

\vspace*{2ex}
\item \emph{Asynchronous Proximal Gradient}:  for $f(x) = \frac{1}{m}\sum_{j=1}^m f^j(x)$ with each $(f^j)$ $\mu$-strongly convex with an $L$-Lipchitz gradient,
        \begin{align}
        \label{davepg}
\tag{DAve-PG}
\left\{ 
\begin{array}{ll}
    u_{k+1} &= \displaystyle \frac{1}{m} \sum_{j=1}^m \left( x_{k-d^j_k} - \gamma \nabla f^j(x_{k-d^j_k}) \right)  \\
    x_{k+1} &= \prox_{\gamma g}(u_{k+1})
\end{array}
\right.
\end{align}
where an iteration is performed as soon as an oracle (say $i$) finishes computing a gradient ($\nabla f^i(x_{k-d_k^i})$); then, $u_{k+1}, x_{k+1}$ are updated. The time $k-d^j_k$ is then the iteration of the last point for which $j$ has sent a response. Provided that no oracle completely stops working, DAve-PG {satisfies} \cref{asm:cv} for $\gamma\in(0,2/(\mu+L)]$ {\cite[Th.~3.2]{mishchenko2018distributed}}.
\end{itemize}

This inventory, far from complete, gives an idea of the diversity of favorable proximal methods. We end this section by mentioning methods that \emph{do not satisfy} \cref{asm:cv}; they mostly fall into two categories: 
\vspace*{-0ex}
\begin{itemize}
    \item \emph{Stochastic versions of proximal gradient.} 
    There are situations where the gradient is obtained via a stochastic oracle $g_k$ {satisfying} $\mathbb{E}[g_k] = \nabla f(x_k)$ and $ \mathbb{E}[ \| g_k - \nabla f(x_k) \|^2]\leq \sigma^2$. The proximal stochastic gradient, obtained by replacing the gradient in \eqref{pg} by such stochastic one,
    {may not recover the structure of the problem; see \cite[Sec.~1.3]{poon2018local}. This case indeed
    fails to satisfy \cref{asm:cv}: there is convergence only in probability and not almost surely (or any other types of weaker convergence giving localization guarantees).}
    
    \vspace*{2ex}
    \item \emph{Inexact proximal methods.} When the proximal operator is not cheap, the step $x_{k+1} = \prox_{\gamma g}(u_{k+1})$ may be computed only approximately (see e.g.\;\cite{rockafellar1976monotone}). The scheme \eqref{algo} is not exactly {satisfied} and thus fall out of the scope of the present paper. In addition, inexact methods may not recover the structure of the optimization problem. 
\end{itemize}

\section{Proximal Identification}
\label{sec:ident}

Proximal methods are the natural algorithms for solving regularized learning problems. The convergence analysis of these algorithms is well-known and studied in the references given in the previous section. Less-known and studied is the qualitative behavior (beyond convergence) of these methods. In this section, we question where the iterates of proximal methods are located and show that they reach some near-optimal structure after some finite time: this is the \emph{identification} property of proximal methods.

\subsection{A generic identification result}

The observation that the iterates produced by optimization algorithms reach some structure dates {can be traced back to e.g.\;\cite{wright1993identifiable} for constrained problems and \cite{hare2004identifying} for nonsmooth minimization.}
Our presentation here is original regarding its simplicity and generality.
For additional details on identification, including early references, we point to\;\cite{fadili2018sensitivity} and references therein. 

In order to properly study structure membership, let us define the \emph{sparsity vector} that describes the structure of a point $x\in\mathbb{R}^n$ relatively to the  collection $\collong$ introduced in the framework \eqref{fw}, as the vector  ${\Str}(x)\in \{0,1\}^\nM$ defined as 
\begin{align}\label{eq:sparse}
    \left[ {\Str}(x) \right]_i = 0 ~~\text{ if } x \in \M_i \text{ and } 1 \text{ elsewhere}.
\end{align}

An identification theorem is a result describing the eventual structure of the iterates of an algorithm, i.e. the value of ${\Str}(x_k)$ for large enough $k$, with respect to the optimal point $x^\star$. We provide below an elementary but powerful enough identification theorem.

\begin{theorem}[Enlarged identification]
\label{th:identification}
Let \cref{asm:cv} hold. Then, the iterates $(x_k)$ produced by \eqref{algo} identify some manifolds in the  collection of the framework \eqref{fw}\footnote{{%
For simplicity, we place ourselves within the framework\;\eqref{fw}, because it is needed to exploit identification in practice. However, the result uses only existence and uniqueness of the proximal mapping (first part of item 3 in the definition of framework \eqref{fw}).}}. More precisely, {for each $\varepsilon>0$ there is a $K$ such that for all $k > K$ we have,}
with probability one,
{for all $i$} %
\begin{equation}\label{eq:identification_result}
	 {[{\Str}(x^\star)]_i ~\leq~ [{\Str}(x_k)]_i ~\leq~ \max\left\{ [{\Str}\big(\mathbf{prox}_{\gamma g}(u)\big)]_i :  \|u-u^\star\|\leq \varepsilon \right\}}%
\end{equation}
where $x^\star$ is the minimizer of \eqref{eq:main_problem} to which $(x_k)$ converges and $u^\star$ is defined as in \cref{asm:cv}.
\end{theorem}

Before going to the proof, we note that a slightly more precise result is shown there: the right part $[{\Str}(x_k)]_i\leq \max\left\{ [{\Str}\big(\mathbf{prox}_{\gamma g}(u)\big)]_i :  \|u-u^\star\|\leq \varepsilon \right\}$ actually holds as soon as $\|u_k-u^\star\|\leq \varepsilon$ .

\begin{proof}
For the right part of the result, since $u_k \to u^\star$ almost surely from \cref{asm:cv}, we have that for any $\varepsilon>0$, $\|u_k-u^\star\|\leq \varepsilon$ for $k$ large enough with probability one.
This directly gives the right hand side of the inequality. Then, ${\Str}\big(\mathbf{prox}_{\gamma g}(u_k)\big) \leq \max\left\{ {\Str}\big(\mathbf{prox}_{\gamma g}(u)\big) :  \| u - u^\star \|\leq \varepsilon \right\}$.

The left part of the inequality is more interesting.  Consider the collection of sets to which $x^\star$ belongs $\col^\star = \{ \M_i\in \col : x^\star \in \M_i \}$ and its complementary collection $\overline{\col^\star} = \col \setminus \col^\star$. Since $\M^\circ := \bigcup \{\M \in \overline{\col^\star}\} $ is a closed set (given that it is a finite\footnote{This argument is the main reason why finiteness is assumed in the framework \eqref{fw}. For countably many sets in the collection, $\M^\circ$ may not be closed (it is then a $F_{\sigma}$ set).} union of closed sets), its complementary set $\mathbb{R}^n\setminus \M^\circ $ is open. Since $x^\star \in \mathbb{R}^n\setminus \M^\circ $ by definition of $\col^\star$, there exists a ball of radius $\varepsilon'>0$ around $x^\star$ that is included in $\mathbb{R}^n\setminus  \M^\circ  $. Hence,  as $x^k \to x^\star $ almost surely, it will reach this ball in finite time with probability one and thus belong to fewer subspaces than $x^\star$. No point $x$ in that ball belongs to $ \M^\circ $ which means that ${\Str}(x^\star) \leq {\Str}(x)$ coordinate-wise.
\end{proof}

\begin{figure}[!h]
     \centering
     \includegraphics[width=0.67\textwidth]{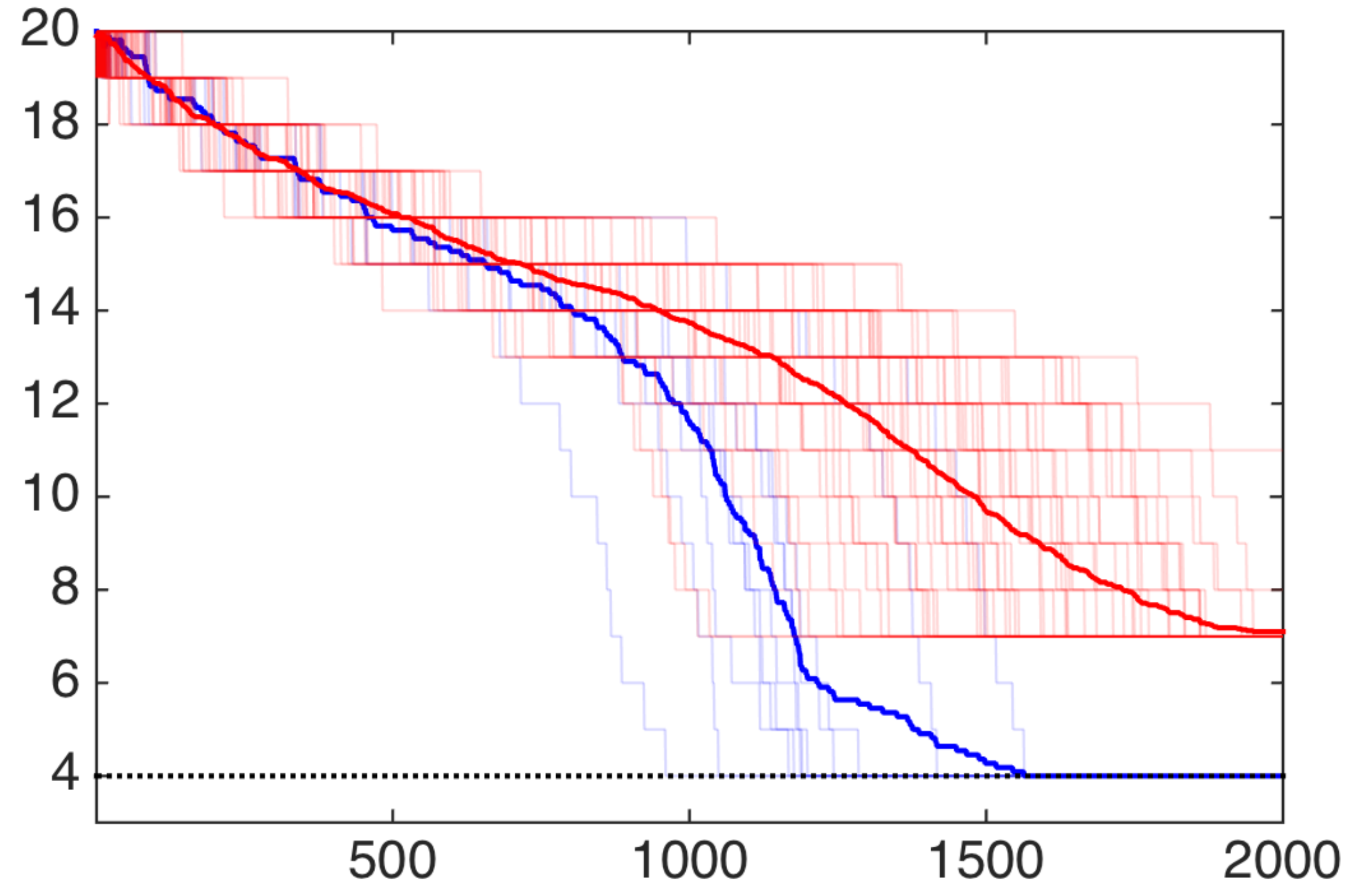}
    \caption{Illustration of the enlarged identification of the proximal gradient algorithm solving matrix least-squares problems regularized by nuclear norm. We generate many random problems for matrices of size $20\times20$ and with optimal solutions of rank $4$. We plot the decrease of the rank of the iterates for two groups of problems: in blue the problems are well-posed and the algorithm identifies the optimal rank; in red the problems are degenerate and the algorithm identifies something bigger. Each curve is a trajectory and the bold curves are the averaged trajectory for the two groups.\label{fig:nuclear}}
\end{figure}

This rather generic result demonstrates that the iterates of proximal methods partially identify some structure as soon as $u_k$ is close enough to $u^\star$.
{After a finite but unknown number of iterations, we have:}
\begin{itemize}
    \item \emph{less structure than $x^\star$}. 
    {The left-hand inequality guarantees that all the identified manifolds contain the optimal manifolds,
    but possibly others (which means that less structure is identified).}
    {Note that this does not depend on the algorithm used.}
    \item \emph{some structure.} The identified structure is not random: it is controlled by the right-hand-side term, which encompasses how the structure changes around the pair\;$(x^\star, u^\star)$. This upper-bound is in general impossible to evaluate a priori (for an exception, see\;\cite{duval2017sparse}). Intuitively, its size captures the difficulty of reaching the structure of\;$x^\star$; see Figure\;\ref{fig:nuclear}.
\end{itemize}

Note that the structure mentioned here is relative to {the solution of \eqref{eq:main_problem} to which the algorithm is converging}. If there are multiple solutions to the problem, they may have different sparsity patterns, but identification still holds. 

\cref{th:identification} {frames the structure of the iterates}; sometimes the terminology \emph{enlarged} or \emph{approximate} identification is used. However, there is no reason for $\Str(x_k)$ to be asymptotically:
stable, monotonous, or close to $\Str(x^\star)$ (or to the right bound),
as evidenced numerically in \cite{fadili2018sensitivity,fadili2019model,bareilles2019interplay}. This may pose some problems in practice if one uses $\Str(x_k)$ as a proxy for the optimal structure. Fortunately, the family of problems for which the identification is \emph{exact}, i.e. $\Str(x_k)  = \Str(x^\star)$, can be characterized, as explained in the next section.

\begin{remark}[Relation with Screening]\label{rem:screening}
Finding near-optimal structures is an important part of machine learning problems. For instance, when $g=\|\cdot\|_1$ (for instance in the lasso problem \cite{tibshirani1996regression}), identification brings information about the nullity of the coordinates of the problem solution, enabling feature selection and/or dimension reduction \cite{tibshirani2012strong}. Discarding null features when solving sparse learning problems is often called \emph{screening}.
With the present notation, this means knowing (an upper-bound on) $\Str(x^\star)$ before convergence in order to discard null features and optimizing on the non-zeros entries. \cref{th:identification} implies that $\Str(x_k)$ is a valid upper-bound, but after some unknown time. 

Observe though that, following the same arguments as in the proof of \cref{th:identification}, we have ${\Str}(x^\star)\leq \max\{{\Str}(\mathbf{prox}_{\gamma g}(u)) : u\in \mathcal{X}\}$ for any set $\mathcal{X} \ni u^\star$. Determining such a \emph{safe region} $\mathcal{X}$ \cite{ghaoui2010safe} (for instance using duality \cite{fercoq2015mind}) and computing the max gives a so-called \emph{safe screening rule}. While this may be computationally intractable in many cases, this proved to be very efficient for solving the the lasso problem and other sparse learning models; for additional details, see e.g. \cite{mairal2010representations,ndiaye2017gap,massias2018celer}.\qed
\end{remark}

\subsection{Exact identification under a qualifying constraint}

The general enlarged identification result of \cref{th:identification} can be easily strengthened to an exact version, guaranteeing $\Str(x_k)  = \Str(x^\star)$ after some time.
A sufficient condition for exact identification is that the proximity operator maps all points $u$ close to $u^\star$ to a point with the same structure as $x^\star$:
\begin{align}
    \tag{QC}
    \label{eq:QC}
   \exists ~ \varepsilon>0 \text{ such that } ~~~ {\Str}(x^\star) =
   \max\left\{{\Str}(\mathbf{prox}_{\gamma g}(u)) : \|u-u^\star\|\leq \varepsilon \right\}
\end{align}

Note that this qualifying condition only depends on the optimal pair $(x^\star,u^\star)$ and the proximity operator of function $g$. 

\begin{corollary}[Exact identification]
\label{cor:identification_exact}
Let \cref{asm:cv} hold; suppose that \eqref{eq:main_problem} has a unique minimizer\footnote{In \cref{th:identification}, the result holds for \emph{the} minimizer to which the algorithm converges. Here, for simplicity, we assume furthermore uniqueness of the minimizer. Nevertheless, if there are multiple minimizers and condition\;\eqref{eq:QC} holds for all pairs $(x^\star,u^\star)$, the structure of all optimal points can be proved to be the same under mild conditions. For instance, if $f$ is differentiable and $\nabla f$ is $L$-Lipchitz continuous on the set of optimal solutions $X^\star$, 
then at optimality $x^\star = \prox_{\gamma g}(u^\star) \text{ with } u^\star = x^\star - \gamma \nabla f(x^\star)$. If \eqref{eq:QC}  holds for some $\varepsilon>0$, take two solutions $x^\star_1,x^\star_2\in X^\star$ such that $\| x^\star_1 - x^\star_2 \| \leq \varepsilon/(1+\gamma L)$; then, $\|u_1^\star - u_2^\star \| \leq \varepsilon$ and thus $\Str(x_1^\star) = \Str(x_2^\star)$ by \eqref{eq:QC}. Then, since $\Str(x^\star)$ is the same for all minimizers, \cref{cor:identification_exact} holds. } and that \eqref{eq:QC} holds true. Then, the iterates $(x_k)$ produced by \eqref{algo} identify the optimal structure in the framework \eqref{fw}. More precisely, for $k$ large enough, we have
\begin{equation}\label{eq:identification_result}
{\Str}(x_k) = {\Str}(x^\star)   ~~~~ \text{ with probability } 1
\end{equation}
where $x^\star$ is the minimizer of \eqref{eq:main_problem} and $u^\star$ is defined as in \cref{asm:cv}.
\end{corollary}

\begin{remark}[Identification time]
As mentioned previously, identification time is generally unknown. However, if \eqref{eq:QC} holds for some known $\varepsilon>0$, then any $u_k$ such $\|u_k-u^\star\|\leq \varepsilon$ will lead to a $x_k$ {satisfying}  $\Str(x_k)  = \Str(x^\star)$. Using the convergence rate of the algorithm, an estimate of the identification time can be given. For instance, if a proximal method {satisfies} $\|u_k - u^\star\| = \mathcal{O}(1/k)$, then we have $\Str(x_k)  = \Str(x^\star)$ after $\mathcal{O}(1/\varepsilon)$ iterations. This rationale is used in\;\cite{garrigos2017thresholding,nutini2019active,sun2019we}.\qed
\end{remark}

\begin{remark}[Link between the qualifying constraint and partial smoothness] Condition \eqref{eq:QC} depends only weakly on the optimal structure, which is in contrast with part of identification literature that considers $g$ partly smooth \cite{hare2004identifying,liang2017activity} relative to the final manifold. As shown in \cite[Apx.\;A]{bareilles2019interplay}, if $g$ is partly smooth relative to the manifold $ \M^\star := \bigcap \{ \M_i\in \col: x^\star\in\M_i \} $ and the condition $(u^\star - x^\star)/\gamma \in \mathrm{ri}~\partial g(x^\star)$ is {satisfied}, then \eqref{eq:QC} holds.
Thus the 
assumption is a consequence of %
the ``irrepresentable condition'' which is a classical assumption in machine learning; see e.g.\;\cite{zhao2006model} and\;\cite{bach2008trace}.
Note finally, that this condition simplifies to $-\nabla f(x^\star) \in \mathrm{ri}~ \partial g(x^\star)$ for the proximal gradient method \eqref{pg} (and most derivatives such as \eqref{apg}, \eqref{davepg}).\qed
\end{remark}

\section{Nonsmoothness can help computationally}
\label{sec:harness}

Nonsmooth regularizations are used in machine learning and signal processing for the nice \emph{recovery} or \emph{consistency} properties\footnote{These properties are intrinsic consequences of the nonsmoothness of regularizations, and directly connected to the stability properties of optimal solutions and identification properties of algorithms. Among a rich literature, we refer to the statistical properties of lasso \cite{tibshirani1996regression}, the recovery in compress sensing literature \cite{candes2006robust,donoho2006compressed}, and the general  model consistency result of~\cite{fadili2019model}. } that they induce; see e.g.\;\cite{vaiter2015low} and recall Figure\;\ref{fig:lasso}. The point of this paper is to advocate that nonsmoothness can also be useful to improve the optimization algorithms used in these applications. This is what we illustrate in this section, building on the tools presented so far.

In the previous section, we have properly formalized that proximal methods uncover (some of) the optimal structure of the optimization problem, which corresponds to the priors put in the learning model. Another important thing to keep in mind is that while we are able to observe the structure along the way, we do not know in general if this structure is in fact optimal or if it will remain stable.
The key question is then
\begin{center}
    \itshape How can we harness the identified structure ?
\end{center}

We present below general ideas to exploit the structure brought by proximal identification.
These ideas can be seen as good practices or, more generally, important points to keep in mind when designing algorithms on problems with noticeable structure.

\subsection{Compression: using the structure to compress iterates}

The first way to harness structure is very simple and direct: since the algorithm's iterates will eventually reach some structure, we may use it to encode them more efficiently. This idea is mainly useful with sparse structures (see the examples of Section\;\ref{sec:NS}): the iterates can be stored in a \texttt{(key,value)} form which directly saves space when storing the iterates. 

\begin{figure}[!h]
     \centering
     \begin{subfigure}[t]{0.49\textwidth}
         \centering
         \resizebox{\textwidth}{!}{         \input{identif/density.tex}      }
         \caption{Iterates density along a run of DAve-PG.}
     \end{subfigure}\hfill
     \begin{subfigure}[t]{0.49\textwidth}
        \centering
        \resizebox{\textwidth}{!}{     \input{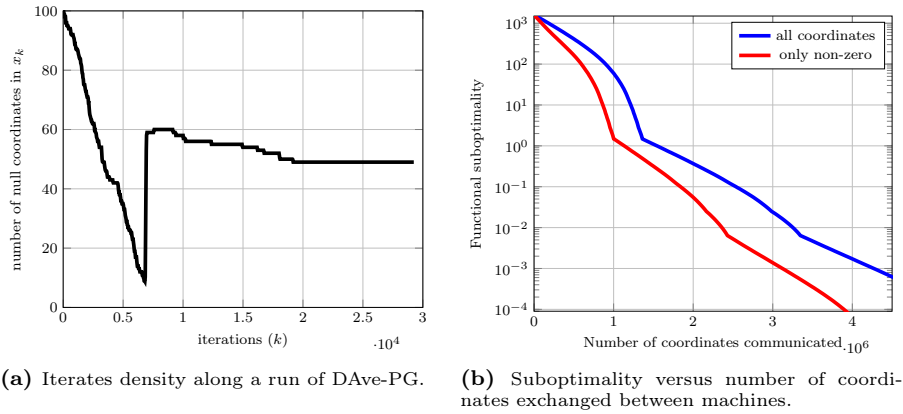}    }
         \caption{Suboptimality versus number of coordinates exchanged between machines.}%
    \end{subfigure}
    \caption{Identification of null coordinates and its use for reducing the exchanges for \eqref{davepg} on a $200\times100$ lasso problem ran on one machine/master coordinating $10$ oracles/workers. We plot on the subfigure (b) the decrease of objective function with respect, not to iterations as usual, but to number of coordinates exchanged between the master and the workers. We see exchanging only non-zero coordinates is directly beneficial in terms of communication cost. The three regimes in the curves correspond to the three behaviors shown on subfigure (a): first the algorithm gives priority to sparsity against the quadratic term, until around iteration 7000, where the huge change in sparsity depicts an approximate identification; and finally exact identification takes place after 2e4 iterations. Recall from Section\;\ref{sec:algo} that an iteration here corresponds to an update from one of the 10 workers. In this distributed framework, the bottleneck is the communications; we show that identification would automatically reduce communications.
    Finally, we mention that this example also shows that identification is not monotone in general: contrary to Figure\;\ref{fig:nuclear}, subfigure (a) shows huge variations in the current support. \label{fig:sparse}}
\end{figure}

More interestingly, in the case of distributed algorithms involving communications between a master and some workers, such as \eqref{davepg}, a sparse encoding directly reduces the communication cost, which is often a bottleneck in distributed systems; see Figure~\ref{fig:sparse} and\;\cite{mishchenko2018delay,grishchenko2018asynchronous}.

\vspace*{-2ex}
\subsection{Acceleration: promoting structure stability over speed}

\begin{figure}[!h]
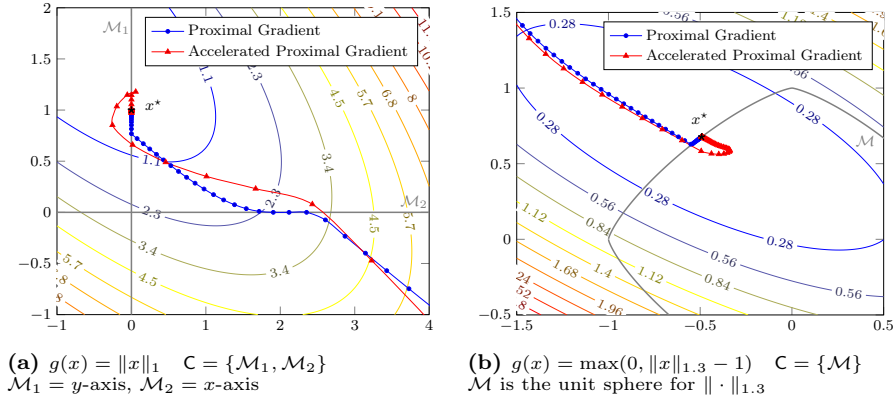

     \centering
     \begin{subfigure}[t]{0.49\textwidth}
         \centering
         \resizebox{\textwidth}{!}{         \input{identif/Lasso-2d-IstaFista_iterates.tikz}      }
         \caption{$g(x) = \|x\|_1 ~~~ \col=\{\M_1,\M_2\}$ \\ ~~~~$\M_1 = y$-axis, $\M_2 = x$-axis}
     \end{subfigure}\hfill
     \begin{subfigure}[t]{0.49\textwidth}
        \centering
        \resizebox{\textwidth}{!}{     \input{identif/distunitball-norm1.3_iterates.tikz}   }
         \caption{ $g(x) =  \max(0, \|x\|_{1.3}-1) ~~~ \col=\{\M\}$ \\ ~~~~$\M$ is the unit sphere for $\|\cdot\|_{1.3}$}
    \end{subfigure}
    \caption{Iterates behavior for the Proximal Gradient (with and without acceleration) when minimizing a function of the type $\|Ax-b\|^2 + g(x)$ with two different functions~$g$, the $\ell_1$-norm (left) and  
    the distance to the unit ball of the 1.3 norm (right). \label{fig:identif}}
\end{figure}

The proximal gradient \eqref{pg} can be accelerated using inertia (see~\eqref{apg}); this acceleration makes the worst case functional convergence rate go from $\mathcal{O}(1/k)$ to $\mathcal{O}(1/k^2)$, which is the optimal complexity for such problems \cite{nesterov1983method}. However, acceleration interferes with identification: we illustrate this graphically on Figure~\ref{fig:identif} which
displays the iterates of \eqref{pg} and \eqref{apg} on two problems~\eqref{eq:main_problem} 
in~$\RR^2$ with $f$ quadratic and $g$ nonsmooth.  
This figure highlights two properties:
\begin{itemize}
    \item Acceleration is faster and has some kind of exploratory behavior which increases its ability to reach the optimal structure (1 iteration = 1 mark in Figure~\ref{fig:identif}). 
    \item Even though the accelerated variant often reaches the optimal structure faster, the additional inertial term in \eqref{apg} is proportional to the difference between the last two iterates. Thus, when reaching a new manifold, the inertial term may drive the iterates out of it (In part (a), \eqref{apg} arrives to the optimal manifold but then is dragged leftwards). More generally, even if \eqref{apg} {satisfies} \cref{asm:cv} and eventually identifies, it is less stable than the vanilla proximal gradient.
\end{itemize}

To summarize, accelerated and vanilla proximal gradient both have respective merits in terms of identification. In order to try and combine these merits, an intermediate method, introduced in 
\cite{bareilles2019interplay},
consists in accelerating only if %
the structure is preserved. This strategy shows encouraging results by being reasonably fast while preserving a good structure identification. This makes it interesting in machine learning problems where recovering (partial) structure is almost as important as convergence.

\subsection{Dimension reduction: updating {with respect to the} identified structure}
\label{sec:alg_ident}

Let us make a simple remark. Once a manifold has been identified (e.g.~some coordinates are null), there is some hope that it will stay identified for the remainder of the algorithm (e.g.~the coordinates will remain null for all subsequent iterations). It is thus natural to update preferentially outside of this manifold (e.g.~update preferentially the non-null coordinates). Mathematically, this leads to the following conceptual algorithm:
\begin{itemize}
    \item Observe $\Str(x_k)$, which gives $\col_k = \{ \M_i\in \col : x_k \in \M_i \}$, the identified collection, and $\M_k = \bigcap_{\M\in\col_k}\M$, the identified manifold. 
    \item Compute a \emph{working manifold} $\N_k$ from $\col_k$ and $\M_k$.
    \item Update iterates preferentially\footnote{Preferentially but not only. Since the user never knows if the optimal structure is attained, the directions outside the identified manifold $\M_k $ should still be updated at times.%
    } along $\N_k$.
\end{itemize}
Two types of strategies for choosing $\N_k$ and updates are presented in the next two subsections: (i) alternating between updates along $\N_k$ and updates on the whole space; (ii) drawing randomly $\N_k$ as a {larger} manifold %
so that, in expectation, the full space is covered.

\subsubsection{Predictor-corrector methods}\label{sec:riemann}

After the observation of the current collection $\col_k$, a simple and natural choice of working manifold is to consider
\begin{align}
    \label{eq:randomup}
    \N_k = %
    \bigcap_{\M\in\col_k} \M \;.
\end{align}
In the ideal situation where we know that our working manifold $\N_k$ is optimal (i.e.\;$x^\star\in\N_k$), we would just have to solve our optimization problem restricted to\;$\N_k$. However, as already mentioned, we never know if $\N_k$ is optimal or not. An idea is therefore to interlace efficient steps along $\N_k$ (as if we knew the optimal manifold) and proximal steps in the whole space (to keep on identifying).
This type of methods is called predictor-corrector in \cite{daniilidis2006geometrical}.

Interestingly, the nonsmooth function $g$ restricted to the manifold $\N_k$ is often smooth, and the problem restricted to\;$\N_k$ is then a Riemannian optimization problem\;\cite{absil-etla-2008}. In words: on top of reducing dimension, identification may also allow us to get rid of nonsmoothness; at the price of the additional %
constraint of lying in\;$\N_k$. See\;\cite{miller2005newton} for discussions on the connections between nonlinear programming, nonsmooth optimization, and Riemannian optimization. 

We write the generic predictor-corrector method exploiting this identification of smooth substructure, as follows
\begin{align}
    \left\{ 
        \begin{array}{rll}
            u_k = & \text{(Gradient or Newton) step %
            tangent to $\N_k$} & ~~~\text{(Riemannian step)} \\[0.2cm]
            x_{k+1} =& \prox_{\gamma g}(u_k) & ~~~\text{(prox)}
        \end{array}
    \right.
\end{align}
When a Riemannian Newton step is used, this method
offers a geometrical interpretation of the so-called
$\mathcal{VU}$ algorithms (introduced in \cite{lemarechal-oustry-sagastizabal-2000} and showed to be superlinearly convergent in \cite{mifflin2005algorithm}).

\subsubsection{Randomized structured descent}

We focus here on~\eqref{eq:main_problem} for $f$ convex and smooth with a Lipschitz gradient, $g$ convex lower semi-continuous, and 
linear  manifolds, as for instance for sparse (Section~\ref{sec:sparsecol}) or flat (Section~\ref{sec:colTV}) structures. 

After the observation of the current collection $\col_k$, \cite{grishchenko2020proximal} proposes to update {on the random working subspace}
\begin{align}
    \N_k = %
    \bigcap_{\M\in\col_k} \left( \xi_{\M,k} \M + (1-\xi_{\M,k}) \mathbb{R}^n \right)
\end{align}
where the $(\xi_{\M,k})$ are independent identically distributed Bernoulli random variables with success probability $p\in(0,1)$. This means that the working subspace is obtained from the current collection, by randomly removing some identified spaces. {This allows any direction in $\RR^n$ to be explored with non-zero probability, which 
prevents the algorithm  
against mis-identification.} 
For instance, in the case of sparsity manifolds, it consists in working on variables with the same support as the current variable, except for some randomly drawn coordinates that we allow to be non-zero as well.

{Working in the subspace $\N_k$ involves projecting the gradient,  %
 which unfortunately adds a bias in the averaged decent direction. To compensate, \cite{grishchenko2020proximal} proposes to debias the iterates with the inverse squared root of the averaged projection $ Q_{k} = \left( \mathbb{E} ~ \proj_{\N_k} \right)^{-\frac{1}{2}}$, which leads to the following algorithm.}

\begin{align}
    \left\{ 
        \begin{array}{rll}
            y_k = & x_k - \gamma \nabla f(x_k) & \text{(gradient)} \\[0.2cm]
            u_k = & Q_k^{-1} \left( \proj_{\N_k}\left( Q_k y_k  \right) +  \proj_{\N_k}^\perp\left( u_{k-1} \right) \right)  & \text{(debiased projection)} \\[0.2cm]
            x_{k+1} =& \prox_{\gamma g}(u_k) & \text{(prox)}   \\[0.2cm]
        \end{array}
    \right.
\end{align}

{
Unfortunately, as such, this method does not work if the collection $\col_k$ (or equivalently the distribution of subspaces $\N_k$) changes too much%
, which goes against the progressive structure identification. To overcome this, \cite{grishchenko2020proximal} proposes to change the current collection $\col_k$ (onto which~$\N_k$ is based) not at every iteration but after some waiting time depending on the amount of change in the iterates' structure (the more changes, the longer the wait).
Thus, by adapting the collection at the right frequency, the method identifies the optimal structure %
 and benefits 
from a competitive rate of convergence compared to the vanilla proximal gradient while using a smaller part of the gradient at each iteration; see \cite{grishchenko2020proximal} for details and illustrations. %
}

\section{Conclusions}

In this paper, we introduced a framework to formalize nonsmooth optimization problems with noticeable structure. These problems often arise in machine learning applications and are prone to be solved by proximal methods. These algorithms are known to identify the problem structure, but going one step further, we show that the structure progressively uncovered by these methods can be harnessed in practice to improve their performances.
Among many possible applications, we emphasize the automatic compression, the interplay with acceleration, and some algorithmic options offered by dimension reduction.

\section*{Acknowledgements}

The authors would like to warmly thank the whole \emph{DAO team} and especially our PhD students Gilles Bareilles, Mathias Chastan, Sélim Chraibi, Dmitry Grishchenko, Yu-Guan Hsieh, and Yassine Laguel. FI benefited from the support of the ANR JCJC project \emph{STROLL} (ANR-19-CE23-0008). This work has been partially supported by MIAI@Grenoble Alpes (ANR-19-P3IA-0003).

\bibliographystyle{spmpsci}
\bibliography{references}

\end{document}

%% file: setup.tex
\usepackage{amsmath}		%
\usepackage{amssymb}		%
\usepackage{amsfonts}		%
\usepackage{amsthm}		    %

\usepackage{mathtools}		%

\usepackage[utf8]{inputenc}		%
\usepackage[T1]{fontenc}		%
\renewcommand\bf{\bfseries}
\renewcommand\tt{\ttfamily}
\renewcommand\it{\itshape}

\usepackage{dsfont}		%

\usepackage{inconsolata}

\usepackage[%
cal=cm,
]
{mathalfa}

\usepackage[labelfont={bf,small},labelsep=colon,font=small, singlelinecheck=false]{caption}	%
\usepackage[dvipsnames,svgnames]{xcolor}		%
\colorlet{MyBlue}{DodgerBlue!75!Black}
\colorlet{MyGreen}{DarkGreen!85!Black}

\usepackage{placeins}

\usepackage{wasysym}		%

\usepackage{fancybox}		%
\usepackage{subcaption}		%
\usepackage{tikz}		%
\usetikzlibrary{calc,patterns}
\usepackage{pgfplots}
\pgfplotsset{compat=1.14}

\usepackage{array}		%
\newcolumntype{C}[1]{>{\centering\let\newline\\\arraybackslash\hspace{0pt}}m{#1}}
\usepackage{booktabs}		%
\usepackage[inline,shortlabels]{enumitem}		%
\usepackage{multirow}
\usepackage{threeparttable}   %

\usepackage{acronym}		%
\usepackage{latexsym}		%
\usepackage{xfrac}		%
\usepackage{xspace}		%
\usepackage{xparse}     %

\usepackage{hyperref}
\hypersetup{
colorlinks=true,
linktocpage=true,
pdfstartview=FitH,
breaklinks=true,
pdfpagemode=UseNone,
pageanchor=true,
pdfpagemode=UseOutlines,
plainpages=false,
bookmarksnumbered,
bookmarksopen=false,
bookmarksopenlevel=1,
hypertexnames=true,
pdfhighlight=/O,
urlcolor=Maroon,linkcolor=MyBlue!60!black,citecolor=MyBlue!70!black,	%
pdftitle={},
pdfauthor={},
pdfsubject={},
pdfkeywords={},
pdfcreator={pdfLaTeX},
pdfproducer={LaTeX with hyperref}
}
\usepackage[sort&compress,capitalize,nameinlink]{cleveref}		%
\crefname{assumption}{Assumption}{Assumptions}
\crefname{assumptionloc}{Assumption}{Assumptions}

\usepackage{autonum}		%

\usepackage[textwidth=20mm]{todonotes}		%
\newtheorem*{corollary*}{Corollary}		%

\newtheorem{assumption}{Assumption}		%

\newtheorem*{definition*}{Definition}		%
\newtheorem*{assumption*}{Assumptions}		%

\makeatletter		%
\newcommand{\asmtag}[1]{%
  \let\oldtheassumption\theassumption%
  \renewcommand{\theassumption}{#1}%
  \g@addto@macro\endassumption{%
    \addtocounter{assumption}{-1}%
    \global\let\theassumption\oldtheassumption}%
  }
\makeatother

\newtheorem{fw}{Framework}

\makeatletter		%
\newcommand{\fwtag}[1]{%
  \let\oldthefw\thefw%
  \renewcommand{\thefw}{#1}%
  \g@addto@macro\endfw{%
    \addtocounter{fw}{-1}%
    \global\let\thefw\oldthefw}%
  }
\makeatother

\newcommand{\proj}{\mathrm{proj}}

\newcommand{\diag}{\mathrm{diag}}

\newcommand{\RR}{\mathbb{R}}

\newcommand{\prox}{\mathbf{prox}}
\newcommand{\algo}{\mathsf{Update}}

\newcommand{\M}{\mathcal{M}}
\newcommand{\Mi}{\mathcal{M}_i}
\newcommand{\nM}{p}

\newcommand{\col}{\mathsf{C}}
\newcommand{\collong}{\mathsf{C}=\{\M_1,\ldots,\M_\nM\}}

\newcommand{\N}{\mathcal{W}}

\newcommand{\Str}{\mathsf{S}}

\newcommand{\risk}{R}

%% file: identif/density.tex
\begin{tikzpicture}
\begin{axis}[grid= major ,
				xlabel = {iterations ($k$)} ,
				ylabel = {number of null coordinates in $x_k$} ,
				xmin = 0, xmax = 30000,
				ymin = 0, ymax = 100,
]
    \addplot+[
       y filter/.code={\pgfmathparse{100-#1}\pgfmathresult},
    line width={2}, solid, no marks, color = black]
        coordinates {
            (0.0,0.0)
            (27.0,0.0)
            (54.0,0.6296296296296333)
            (95.0,1.0)
            (125.0,1.36666666666666)
            (156.0,2.0)
            (188.0,2.03125)
            (218.0,4.0)
            (248.0,4.0)
            (286.0,4.0)
            (310.0,4.958333333333329)
            (343.0,5.0)
            (368.0,5.0)
            (393.0,5.0)
            (432.0,5.410256410256409)
            (461.0,6.0)
            (486.0,6.0)
            (511.0,6.0)
            (540.0,6.0)
            (569.0,6.0)
            (594.0,6.0)
            (619.0,6.560000000000002)
            (644.0,7.760000000000005)
            (673.0,8.0)
            (705.0,8.0)
            (731.0,8.0)
            (762.0,8.0)
            (793.0,8.0)
            (817.0,8.0)
            (847.0,8.599999999999994)
            (879.0,9.0)
            (907.0,9.32142857142857)
            (935.0,10.107142857142861)
            (961.0,11.0)
            (988.0,11.0)
            (1011.0,11.695652173913047)
            (1042.0,12.0)
            (1068.0,12.84615384615384)
            (1090.0,13.0)
            (1119.0,13.0)
            (1152.0,13.0)
            (1185.0,13.0)
            (1219.0,13.32352941176471)
            (1250.0,14.0)
            (1282.0,14.75)
            (1308.0,15.0)
            (1339.0,15.709677419354833)
            (1365.0,16.0)
            (1404.0,16.0)
            (1435.0,16.225806451612897)
            (1469.0,17.705882352941174)
            (1500.0,18.0)
            (1533.0,18.727272727272734)
            (1561.0,20.0)
            (1590.0,21.0)
            (1621.0,21.16129032258064)
            (1645.0,22.5)
            (1679.0,23.058823529411768)
            (1714.0,24.828571428571422)
            (1743.0,25.0)
            (1774.0,25.0)
            (1795.0,25.0)
            (1822.0,25.77777777777777)
            (1851.0,26.551724137931032)
            (1877.0,27.0)
            (1915.0,27.078947368421055)
            (1941.0,28.0)
            (1969.0,28.5)
            (1995.0,29.0)
            (2027.0,29.0)
            (2052.0,29.0)
            (2084.0,29.4375)
            (2110.0,31.0)
            (2139.0,31.0)
            (2167.0,32.71428571428572)
            (2198.0,34.2258064516129)
            (2222.0,35.0)
            (2255.0,35.45454545454547)
            (2289.0,36.0)
            (2319.0,36.53333333333333)
            (2349.0,37.0)
            (2374.0,37.0)
            (2406.0,37.5625)
            (2444.0,38.0)
            (2472.0,38.0)
            (2505.0,38.0)
            (2535.0,38.0)
            (2561.0,38.0)
            (2594.0,38.75757575757575)
            (2618.0,40.0)
            (2649.0,40.0)
            (2675.0,40.84615384615384)
            (2714.0,41.0)
            (2749.0,41.31428571428572)
            (2787.0,42.60526315789474)
            (2822.0,43.0)
            (2849.0,43.18518518518519)
            (2876.0,44.0)
            (2920.0,44.0)
            (2955.0,44.0)
            (2988.0,45.121212121212125)
            (3022.0,46.0)
            (3047.0,46.0)
            (3077.0,46.0)
            (3102.0,46.0)
            (3131.0,46.0)
            (3162.0,46.0)
            (3190.0,46.0)
            (3219.0,48.0)
            (3245.0,49.84615384615384)
            (3277.0,50.84375)
            (3310.0,51.0)
            (3340.0,51.0)
            (3363.0,51.0)
            (3384.0,51.0)
            (3418.0,51.029411764705884)
            (3447.0,52.0)
            (3479.0,52.9375)
            (3511.0,53.6875)
            (3548.0,54.0)
            (3574.0,54.0)
            (3605.0,54.0)
            (3634.0,54.17241379310346)
            (3667.0,55.24242424242425)
            (3696.0,56.0)
            (3719.0,56.0)
            (3753.0,56.0)
            (3784.0,56.0)
            (3810.0,56.0)
            (3832.0,56.0)
            (3859.0,56.22222222222223)
            (3893.0,57.0)
            (3924.0,57.0)
            (3950.0,57.0)
            (3979.0,57.0)
            (4010.0,57.0)
            (4045.0,57.0)
            (4072.0,57.0)
            (4098.0,57.0)
            (4129.0,57.0)
            (4155.0,57.69230769230768)
            (4181.0,58.0)
            (4208.0,58.0)
            (4235.0,58.0)
            (4263.0,58.0)
            (4292.0,58.0)
            (4326.0,58.0)
            (4354.0,58.0)
            (4381.0,58.0)
            (4411.0,58.0)
            (4432.0,58.0)
            (4469.0,58.0)
            (4498.0,58.0)
            (4525.0,58.0)
            (4557.0,58.21875)
            (4588.0,59.0)
            (4614.0,60.38461538461539)
            (4638.0,61.0)
            (4662.0,61.0)
            (4690.0,61.96428571428572)
            (4722.0,62.0)
            (4750.0,62.0)
            (4785.0,62.97142857142856)
            (4809.0,64.0)
            (4836.0,64.0)
            (4871.0,64.08571428571429)
            (4905.0,65.0)
            (4932.0,65.0)
            (4963.0,65.0)
            (4992.0,65.0)
            (5021.0,65.44827586206895)
            (5045.0,66.29166666666666)
            (5077.0,67.0)
            (5097.0,67.0)
            (5124.0,67.0)
            (5154.0,68.0)
            (5182.0,69.0)
            (5214.0,69.53125)
            (5241.0,70.03703703703704)
            (5268.0,71.0)
            (5298.0,71.0)
            (5322.0,71.0)
            (5343.0,72.0)
            (5370.0,72.07407407407408)
            (5399.0,73.0)
            (5430.0,73.0)
            (5453.0,73.0)
            (5478.0,73.72)
            (5507.0,74.0)
            (5535.0,74.0)
            (5568.0,74.0)
            (5598.0,74.0)
            (5634.0,74.30555555555554)
            (5663.0,75.0)
            (5684.0,75.0)
            (5711.0,75.0)
            (5738.0,75.96296296296296)
            (5780.0,77.0)
            (5802.0,77.0)
            (5833.0,77.0)
            (5867.0,78.23529411764707)
            (5892.0,79.0)
            (5919.0,80.18518518518519)
            (5943.0,81.0)
            (5974.0,81.0)
            (6001.0,81.55555555555554)
            (6031.0,83.0)
            (6058.0,83.0)
            (6082.0,83.0)
            (6110.0,83.0)
            (6149.0,84.71794871794873)
            (6185.0,86.0)
            (6212.0,86.0)
            (6236.0,86.0)
            (6256.0,86.0)
            (6281.0,86.84)
            (6322.0,87.0)
            (6349.0,87.0)
            (6371.0,87.36363636363637)
            (6400.0,88.0)
            (6434.0,88.0)
            (6467.0,88.0)
            (6490.0,88.0)
            (6514.0,88.0)
            (6543.0,88.0)
            (6572.0,88.41379310344828)
            (6597.0,89.0)
            (6627.0,89.0)
            (6655.0,89.0)
            (6688.0,89.87878787878788)
            (6718.0,90.0)
            (6753.0,90.97142857142856)
            (6784.0,91.0)
            (6811.0,91.18518518518519)
            (6844.0,89.75757575757575)
            (6870.0,77.15384615384616)
            (6902.0,60.71875)
            (6931.0,46.896551724137936)
            (6956.0,41.91999999999999)
            (6994.0,41.0)
            (7022.0,41.0)
            (7060.0,41.0)
            (7084.0,41.0)
            (7109.0,41.0)
            (7136.0,41.0)
            (7168.0,41.0)
            (7198.0,41.0)
            (7232.0,41.0)
            (7265.0,41.0)
            (7290.0,41.0)
            (7316.0,41.0)
            (7345.0,41.0)
            (7380.0,41.0)
            (7410.0,41.0)
            (7436.0,41.0)
            (7461.0,41.0)
            (7495.0,41.0)
            (7528.0,41.0)
            (7565.0,40.351351351351354)
            (7590.0,40.0)
            (7622.0,40.0)
            (7650.0,40.0)
            (7684.0,40.0)
            (7716.0,40.0)
            (7747.0,40.0)
            (7773.0,40.0)
            (7801.0,40.0)
            (7829.0,40.0)
            (7856.0,40.0)
            (7888.0,40.0)
            (7918.0,40.0)
            (7949.0,40.0)
            (7979.0,40.0)
            (8004.0,40.0)
            (8034.0,40.0)
            (8061.0,40.0)
            (8089.0,40.0)
            (8125.0,40.0)
            (8155.0,40.0)
            (8179.0,40.0)
            (8211.0,40.0)
            (8235.0,40.0)
            (8263.0,40.0)
            (8292.0,40.0)
            (8319.0,40.0)
            (8346.0,40.0)
            (8374.0,40.0)
            (8413.0,40.0)
            (8438.0,40.0)
            (8463.0,40.0)
            (8488.0,40.0)
            (8519.0,40.0)
            (8553.0,40.0)
            (8582.0,40.0)
            (8616.0,40.0)
            (8646.0,40.0)
            (8676.0,40.0)
            (8707.0,40.0)
            (8734.0,40.0)
            (8757.0,40.0)
            (8787.0,40.0)
            (8821.0,40.0)
            (8852.0,40.0)
            (8878.0,40.0)
            (8918.0,40.0)
            (8948.0,40.0)
            (8980.0,40.0)
            (9014.0,40.0)
            (9041.0,40.0)
            (9073.0,40.0)
            (9102.0,40.0)
            (9127.0,40.0)
            (9153.0,40.61538461538461)
            (9192.0,41.0)
            (9223.0,41.0)
            (9256.0,41.0)
            (9285.0,41.0)
            (9320.0,41.0)
            (9347.0,41.0)
            (9381.0,41.20588235294119)
            (9402.0,42.0)
            (9430.0,42.0)
            (9461.0,42.0)
            (9488.0,42.0)
            (9525.0,42.0)
            (9550.0,42.0)
            (9573.0,42.0)
            (9605.0,42.0)
            (9645.0,42.0)
            (9674.0,42.0)
            (9699.0,42.0)
            (9723.0,42.0)
            (9755.0,42.0)
            (9791.0,42.0)
            (9820.0,42.0)
            (9850.0,42.0)
            (9881.0,42.0)
            (9909.0,42.0)
            (9941.0,42.0)
            (9967.0,42.0)
            (10000.0,42.0)
            (10030.0,43.0)
            (10059.0,43.0)
            (10086.0,43.0)
            (10113.0,43.0)
            (10148.0,43.0)
            (10187.0,43.17948717948718)
            (10217.0,44.0)
            (10251.0,44.0)
            (10285.0,44.0)
            (10314.0,44.0)
            (10337.0,44.0)
            (10363.0,44.0)
            (10403.0,44.0)
            (10434.0,44.0)
            (10458.0,44.0)
            (10486.0,44.0)
            (10518.0,44.0)
            (10542.0,44.0)
            (10573.0,44.0)
            (10606.0,44.0)
            (10641.0,44.0)
            (10669.0,44.0)
            (10698.0,44.0)
            (10724.0,44.0)
            (10753.0,44.0)
            (10786.0,44.0)
            (10821.0,44.0)
            (10848.0,44.0)
            (10876.0,44.0)
            (10903.0,44.0)
            (10932.0,44.0)
            (10961.0,44.0)
            (10987.0,44.0)
            (11019.0,44.0)
            (11040.0,44.0)
            (11065.0,44.0)
            (11097.0,44.0)
            (11125.0,44.0)
            (11151.0,44.0)
            (11185.0,44.0)
            (11214.0,44.0)
            (11243.0,44.0)
            (11276.0,44.0)
            (11308.0,44.0)
            (11340.0,44.0)
            (11362.0,44.0)
            (11396.0,44.0)
            (11426.0,44.0)
            (11448.0,44.0)
            (11475.0,44.0)
            (11502.0,44.0)
            (11534.0,44.0)
            (11567.0,44.0)
            (11589.0,44.0)
            (11622.0,44.0)
            (11651.0,44.0)
            (11685.0,44.0)
            (11716.0,44.0)
            (11743.0,44.0)
            (11764.0,44.0)
            (11793.0,44.0)
            (11827.0,44.0)
            (11858.0,44.0)
            (11885.0,44.0)
            (11906.0,44.0)
            (11935.0,44.0)
            (11965.0,44.0)
            (11994.0,44.0)
            (12022.0,44.0)
            (12049.0,44.0)
            (12076.0,44.0)
            (12106.0,44.0)
            (12139.0,44.0)
            (12166.0,44.0)
            (12195.0,44.0)
            (12221.0,44.0)
            (12252.0,44.0)
            (12279.0,44.0)
            (12306.0,44.0)
            (12341.0,44.02857142857144)
            (12372.0,45.0)
            (12400.0,45.0)
            (12430.0,45.0)
            (12461.0,45.0)
            (12490.0,45.0)
            (12521.0,45.0)
            (12546.0,45.0)
            (12579.0,45.0)
            (12607.0,45.0)
            (12634.0,45.0)
            (12660.0,45.0)
            (12689.0,45.0)
            (12715.0,45.0)
            (12742.0,45.0)
            (12769.0,45.0)
            (12795.0,45.0)
            (12827.0,45.0)
            (12854.0,45.0)
            (12886.0,45.0)
            (12925.0,45.0)
            (12952.0,45.0)
            (12988.0,45.0)
            (13012.0,45.0)
            (13041.0,45.0)
            (13062.0,45.0)
            (13093.0,45.0)
            (13126.0,45.0)
            (13153.0,45.0)
            (13183.0,45.0)
            (13206.0,45.0)
            (13235.0,45.0)
            (13263.0,45.0)
            (13286.0,45.0)
            (13311.0,45.0)
            (13337.0,45.0)
            (13363.0,45.0)
            (13397.0,45.0)
            (13421.0,45.0)
            (13447.0,45.0)
            (13473.0,45.0)
            (13499.0,45.0)
            (13522.0,45.0)
            (13552.0,45.0)
            (13577.0,45.0)
            (13607.0,45.0)
            (13637.0,45.0)
            (13670.0,45.0)
            (13695.0,45.0)
            (13736.0,45.0)
            (13768.0,45.0)
            (13795.0,45.0)
            (13819.0,45.0)
            (13843.0,45.0)
            (13868.0,45.0)
            (13894.0,45.0)
            (13917.0,45.0)
            (13943.0,45.0)
            (13973.0,45.0)
            (14001.0,45.0)
            (14027.0,45.0)
            (14054.0,45.0)
            (14082.0,45.0)
            (14107.0,45.0)
            (14138.0,45.0)
            (14165.0,45.0)
            (14189.0,45.0)
            (14220.0,45.0)
            (14248.0,45.0)
            (14277.0,45.0)
            (14308.0,45.0)
            (14336.0,45.0)
            (14361.0,45.0)
            (14396.0,45.0)
            (14425.0,45.0)
            (14452.0,45.0)
            (14486.0,45.0)
            (14518.0,45.0)
            (14547.0,45.0)
            (14583.0,45.0)
            (14612.0,45.0)
            (14641.0,45.0)
            (14666.0,45.0)
            (14695.0,45.0)
            (14723.0,45.0)
            (14752.0,45.0)
            (14778.0,45.0)
            (14807.0,45.0)
            (14839.0,45.0)
            (14864.0,45.0)
            (14895.0,45.0)
            (14919.0,45.0)
            (14946.0,45.66666666666666)
            (14978.0,46.0)
            (15007.0,46.0)
            (15046.0,46.0)
            (15073.0,46.0)
            (15101.0,46.0)
            (15131.0,46.0)
            (15159.0,46.0)
            (15190.0,46.0)
            (15215.0,46.0)
            (15242.0,46.0)
            (15267.0,46.0)
            (15291.0,46.0)
            (15317.0,46.0)
            (15344.0,46.0)
            (15377.0,46.0)
            (15405.0,46.0)
            (15430.0,46.0)
            (15453.0,46.0)
            (15479.0,46.0)
            (15506.0,46.0)
            (15533.0,46.0)
            (15567.0,46.0)
            (15600.0,46.0)
            (15632.0,46.0)
            (15655.0,46.0)
            (15684.0,46.0)
            (15711.0,46.0)
            (15744.0,46.0)
            (15775.0,46.0)
            (15804.0,46.0)
            (15826.0,46.0)
            (15856.0,46.0)
            (15888.0,46.0)
            (15911.0,46.0)
            (15946.0,46.0)
            (15976.0,46.0)
            (16002.0,46.0)
            (16027.0,46.0)
            (16050.0,46.0)
            (16075.0,46.0)
            (16110.0,46.0)
            (16137.0,46.0)
            (16172.0,46.0)
            (16197.0,46.47999999999999)
            (16220.0,47.0)
            (16249.0,47.0)
            (16281.0,47.0)
            (16317.0,47.0)
            (16342.0,47.0)
            (16368.0,47.0)
            (16398.0,47.0)
            (16435.0,47.0)
            (16465.0,47.0)
            (16493.0,47.0)
            (16529.0,47.0)
            (16564.0,47.0)
            (16591.0,47.0)
            (16631.0,47.0)
            (16659.0,47.0)
            (16689.0,47.0)
            (16714.0,47.0)
            (16741.0,47.296296296296305)
            (16770.0,48.0)
            (16795.0,48.0)
            (16827.0,48.0)
            (16862.0,48.0)
            (16889.0,48.0)
            (16913.0,48.0)
            (16947.0,48.0)
            (16975.0,48.0)
            (17008.0,48.0)
            (17039.0,48.0)
            (17074.0,48.0)
            (17100.0,48.0)
            (17132.0,48.0)
            (17163.0,48.0)
            (17194.0,48.0)
            (17225.0,48.0)
            (17251.0,48.0)
            (17280.0,48.0)
            (17310.0,48.0)
            (17339.0,48.0)
            (17363.0,48.0)
            (17397.0,48.0)
            (17422.0,48.0)
            (17453.0,48.0)
            (17489.0,48.0)
            (17516.0,48.0)
            (17544.0,48.0)
            (17569.0,48.0)
            (17600.0,48.0)
            (17624.0,48.0)
            (17653.0,48.0)
            (17680.0,48.0)
            (17709.0,48.0)
            (17738.0,48.0)
            (17765.0,48.0)
            (17792.0,48.0)
            (17818.0,48.0)
            (17846.0,48.0)
            (17876.0,48.0)
            (17906.0,48.0)
            (17938.0,48.0)
            (17969.0,48.0)
            (18004.0,48.0)
            (18030.0,48.65384615384616)
            (18062.0,49.65625)
            (18087.0,50.0)
            (18113.0,50.0)
            (18142.0,50.0)
            (18168.0,50.0)
            (18199.0,50.0)
            (18238.0,50.0)
            (18264.0,50.0)
            (18296.0,50.0)
            (18325.0,50.0)
            (18349.0,50.0)
            (18381.0,50.0)
            (18406.0,50.0)
            (18430.0,50.0)
            (18457.0,50.0)
            (18488.0,50.0)
            (18517.0,50.0)
            (18541.0,50.0)
            (18572.0,50.0)
            (18598.0,50.0)
            (18624.0,50.0)
            (18658.0,50.0)
            (18687.0,50.0)
            (18718.0,50.0)
            (18749.0,50.0)
            (18779.0,50.0)
            (18801.0,50.0)
            (18827.0,50.0)
            (18853.0,50.0)
            (18878.0,50.0)
            (18904.0,50.0)
            (18935.0,50.0)
            (18963.0,50.0)
            (18995.0,50.0)
            (19028.0,50.0)
            (19058.0,50.0)
            (19097.0,50.0)
            (19121.0,50.0)
            (19150.0,50.13793103448276)
            (19188.0,51.0)
            (19217.0,51.0)
            (19249.0,51.0)
            (19276.0,51.0)
            (19302.0,51.0)
            (19329.0,51.0)
            (19361.0,51.0)
            (19393.0,51.0)
            (19420.0,51.0)
            (19448.0,51.0)
            (19480.0,51.0)
            (19503.0,51.0)
            (19535.0,51.0)
            (19568.0,51.0)
            (19602.0,51.0)
            (19643.0,51.0)
            (19670.0,51.0)
            (19700.0,51.0)
            (19732.0,51.0)
            (19766.0,51.0)
            (19794.0,51.0)
            (19818.0,51.0)
            (19851.0,51.0)
            (19881.0,51.0)
            (19908.0,51.0)
            (19938.0,51.0)
            (19965.0,51.0)
            (19999.0,51.0)
            (20024.0,51.0)
            (20053.0,51.0)
            (20076.0,51.0)
            (20104.0,51.0)
            (20127.0,51.0)
            (20150.0,51.0)
            (20181.0,51.0)
            (20207.0,51.0)
            (20233.0,51.0)
            (20262.0,51.0)
            (20291.0,51.0)
            (20318.0,51.0)
            (20347.0,51.0)
            (20377.0,51.0)
            (20409.0,51.0)
            (20441.0,51.0)
            (20467.0,51.0)
            (20494.0,51.0)
            (20526.0,51.0)
            (20553.0,51.0)
            (20579.0,51.0)
            (20606.0,51.0)
            (20629.0,51.0)
            (20658.0,51.0)
            (20687.0,51.0)
            (20719.0,51.0)
            (20749.0,51.0)
            (20780.0,51.0)
            (20808.0,51.0)
            (20837.0,51.0)
            (20863.0,51.0)
            (20889.0,51.0)
            (20922.0,51.0)
            (20950.0,51.0)
            (20990.0,51.0)
            (21015.0,51.0)
            (21041.0,51.0)
            (21067.0,51.0)
            (21096.0,51.0)
            (21122.0,51.0)
            (21152.0,51.0)
            (21185.0,51.0)
            (21217.0,51.0)
            (21248.0,51.0)
            (21277.0,51.0)
            (21307.0,51.0)
            (21334.0,51.0)
            (21362.0,51.0)
            (21392.0,51.0)
            (21422.0,51.0)
            (21457.0,51.0)
            (21487.0,51.0)
            (21511.0,51.0)
            (21536.0,51.0)
            (21564.0,51.0)
            (21588.0,51.0)
            (21617.0,51.0)
            (21642.0,51.0)
            (21668.0,51.0)
            (21701.0,51.0)
            (21728.0,51.0)
            (21756.0,51.0)
            (21790.0,51.0)
            (21815.0,51.0)
            (21843.0,51.0)
            (21880.0,51.0)
            (21925.0,51.0)
            (21957.0,51.0)
            (21983.0,51.0)
            (22011.0,51.0)
            (22037.0,51.0)
            (22067.0,51.0)
            (22090.0,51.0)
            (22121.0,51.0)
            (22146.0,51.0)
            (22174.0,51.0)
            (22198.0,51.0)
            (22233.0,51.0)
            (22258.0,51.0)
            (22284.0,51.0)
            (22314.0,51.0)
            (22350.0,51.0)
            (22378.0,51.0)
            (22406.0,51.0)
            (22442.0,51.0)
            (22477.0,51.0)
            (22504.0,51.0)
            (22542.0,51.0)
            (22574.0,51.0)
            (22608.0,51.0)
            (22637.0,51.0)
            (22664.0,51.0)
            (22687.0,51.0)
            (22718.0,51.0)
            (22752.0,51.0)
            (22782.0,51.0)
            (22811.0,51.0)
            (22837.0,51.0)
            (22873.0,51.0)
            (22900.0,51.0)
            (22940.0,51.0)
            (22964.0,51.0)
            (22987.0,51.0)
            (23013.0,51.0)
            (23046.0,51.0)
            (23078.0,51.0)
            (23111.0,51.0)
            (23145.0,51.0)
            (23168.0,51.0)
            (23196.0,51.0)
            (23225.0,51.0)
            (23255.0,51.0)
            (23285.0,51.0)
            (23310.0,51.0)
            (23335.0,51.0)
            (23367.0,51.0)
            (23395.0,51.0)
            (23428.0,51.0)
            (23451.0,51.0)
            (23479.0,51.0)
            (23513.0,51.0)
            (23544.0,51.0)
            (23575.0,51.0)
            (23609.0,51.0)
            (23639.0,51.0)
            (23675.0,51.0)
            (23706.0,51.0)
            (23740.0,51.0)
            (23770.0,51.0)
            (23804.0,51.0)
            (23832.0,51.0)
            (23869.0,51.0)
            (23899.0,51.0)
            (23934.0,51.0)
            (23965.0,51.0)
            (23992.0,51.0)
            (24022.0,51.0)
            (24051.0,51.0)
            (24076.0,51.0)
            (24101.0,51.0)
            (24137.0,51.0)
            (24161.0,51.0)
            (24195.0,51.0)
            (24231.0,51.0)
            (24263.0,51.0)
            (24289.0,51.0)
            (24317.0,51.0)
            (24351.0,51.0)
            (24378.0,51.0)
            (24402.0,51.0)
            (24431.0,51.0)
            (24464.0,51.0)
            (24490.0,51.0)
            (24523.0,51.0)
            (24551.0,51.0)
            (24583.0,51.0)
            (24611.0,51.0)
            (24639.0,51.0)
            (24672.0,51.0)
            (24706.0,51.0)
            (24732.0,51.0)
            (24760.0,51.0)
            (24785.0,51.0)
            (24810.0,51.0)
            (24848.0,51.0)
            (24882.0,51.0)
            (24916.0,51.0)
            (24952.0,51.0)
            (24983.0,51.0)
            (25011.0,51.0)
            (25038.0,51.0)
            (25067.0,51.0)
            (25100.0,51.0)
            (25129.0,51.0)
            (25154.0,51.0)
            (25181.0,51.0)
            (25210.0,51.0)
            (25234.0,51.0)
            (25261.0,51.0)
            (25298.0,51.0)
            (25326.0,51.0)
            (25358.0,51.0)
            (25387.0,51.0)
            (25414.0,51.0)
            (25448.0,51.0)
            (25478.0,51.0)
            (25507.0,51.0)
            (25536.0,51.0)
            (25561.0,51.0)
            (25593.0,51.0)
            (25626.0,51.0)
            (25655.0,51.0)
            (25684.0,51.0)
            (25709.0,51.0)
            (25741.0,51.0)
            (25773.0,51.0)
            (25804.0,51.0)
            (25834.0,51.0)
            (25861.0,51.0)
            (25901.0,51.0)
            (25932.0,51.0)
            (25958.0,51.0)
            (25986.0,51.0)
            (26014.0,51.0)
            (26039.0,51.0)
            (26071.0,51.0)
            (26101.0,51.0)
            (26131.0,51.0)
            (26156.0,51.0)
            (26181.0,51.0)
            (26215.0,51.0)
            (26240.0,51.0)
            (26269.0,51.0)
            (26305.0,51.0)
            (26337.0,51.0)
            (26361.0,51.0)
            (26400.0,51.0)
            (26429.0,51.0)
            (26462.0,51.0)
            (26490.0,51.0)
            (26513.0,51.0)
            (26545.0,51.0)
            (26576.0,51.0)
            (26610.0,51.0)
            (26639.0,51.0)
            (26666.0,51.0)
            (26696.0,51.0)
            (26730.0,51.0)
            (26758.0,51.0)
            (26787.0,51.0)
            (26811.0,51.0)
            (26838.0,51.0)
            (26867.0,51.0)
            (26899.0,51.0)
            (26927.0,51.0)
            (26960.0,51.0)
            (26986.0,51.0)
            (27017.0,51.0)
            (27047.0,51.0)
            (27079.0,51.0)
            (27116.0,51.0)
            (27145.0,51.0)
            (27177.0,51.0)
            (27203.0,51.0)
            (27234.0,51.0)
            (27264.0,51.0)
            (27293.0,51.0)
            (27317.0,51.0)
            (27341.0,51.0)
            (27371.0,51.0)
            (27404.0,51.0)
            (27433.0,51.0)
            (27459.0,51.0)
            (27489.0,51.0)
            (27524.0,51.0)
            (27560.0,51.0)
            (27590.0,51.0)
            (27619.0,51.0)
            (27651.0,51.0)
            (27682.0,51.0)
            (27711.0,51.0)
            (27736.0,51.0)
            (27766.0,51.0)
            (27800.0,51.0)
            (27827.0,51.0)
            (27852.0,51.0)
            (27874.0,51.0)
            (27904.0,51.0)
            (27939.0,51.0)
            (27965.0,51.0)
            (27993.0,51.0)
            (28029.0,51.0)
            (28061.0,51.0)
            (28086.0,51.0)
            (28110.0,51.0)
            (28144.0,51.0)
            (28166.0,51.0)
            (28190.0,51.0)
            (28220.0,51.0)
            (28245.0,51.0)
            (28277.0,51.0)
            (28304.0,51.0)
            (28334.0,51.0)
            (28362.0,51.0)
            (28391.0,51.0)
            (28418.0,51.0)
            (28448.0,51.0)
            (28478.0,51.0)
            (28504.0,51.0)
            (28533.0,51.0)
            (28565.0,51.0)
            (28590.0,51.0)
            (28617.0,51.0)
            (28642.0,51.0)
            (28671.0,51.0)
            (28697.0,51.0)
            (28729.0,51.0)
            (28757.0,51.0)
            (28786.0,51.0)
            (28823.0,51.0)
            (28852.0,51.0)
            (28877.0,51.0)
            (28909.0,51.0)
            (28939.0,51.0)
            (28967.0,51.0)
            (28988.0,51.0)
            (29018.0,51.0)
            (29042.0,51.0)
            (29078.0,51.0)
            (29102.0,51.0)
            (29123.0,51.0)
            (29149.0,51.0)
            (29184.0,51.0)
            (29210.0,51.0)
            (29243.0,51.0)
        }
        ;
    \addlegendentry{Dave-PG}
    
    \legend{}; %
\end{axis}
\end{tikzpicture}